\documentclass[12pt]{article}

\usepackage{amssymb}
\usepackage{amscd}
\usepackage{epsfig}
\usepackage{graphicx}
\usepackage{theorem}
\usepackage{amsmath}

\newtheorem{thm}{Theorem}[section]

\newtheorem{prop}[thm]{Proposition}
\newtheorem{Theo}{Theorem}

\newtheorem{coro}[thm]{Corollary}

\newtheorem{lemm}[thm]{Lemma}

\newtheorem{assu}[thm]{Assumption}

{\theorembodyfont{\rm}

\newtheorem{defn}[thm]{Definition}
\newtheorem{exam}[thm]{Example}
\newtheorem{rem}[thm]{Remark}
\newtheorem{nota}[thm]{Notation}
}

\newenvironment{proof}{\noindent {\em Proof.}}{\hfill 
$\square$\newline}

\newcommand{\no}{\noindent}

\newcommand{\La}{\Lambda}
\newcommand{\la}{\lambda}
\newcommand{\al}{\alpha}
\newcommand{\C}{\mathbb C}
\newcommand{\Q}{\mathbb Q}
\newcommand{\Z}{\mathbb Z}
\newcommand{\R}{\mathbb R}

\def\cX{{\mathcal X}}
\def\cA{{\mathcal A}}
\def\cB{{\mathcal B}}
\def\cC{{\mathcal C}}

\def\cH{{\mathcal H}}

\def\cK{{\mathcal K}}

\def\cL{{\mathcal L}}

\def\cH{{\mathcal H}}

\def\cS{{\mathcal S}}

\def\bF{{\bf F}}

\def\sT{{\mathsf T}}
\def\sL{{\mathsf L}}
\def\sC{{\mathsf C}}
\def\sc{{\mathsf c}}

\def\sn{{\mathsf n}}
\def\sm{{\mathsf m}}

\def\sK{{\mathsf K}}

\def\ga{{\mathbf G}_a}

\def\dv{{\rm div}}

\newcommand{\G}{\mathrm G}

\newcommand{\PGL}{{\rm PGL}} 
\def\Val{{\rm Val}}

\def\Spec{{\rm Spec}}

\def\ra{\rightarrow}

\def\A{{\mathbb A}}
\def\C{{\mathbb C}}

\def\P{{\mathbb P}}
\def\Q{{\mathbb Q}}

\def\Z{{\mathbb Z}}
\def\C{{\mathbb C}}

\def\bG{{\mathrm G}}

\def\bH{{\mathrm H}}
\def\bH{{\mathrm H}}
\def\bU{{\mathrm U}}
\def\bW{{\mathrm W}}
\def\bZ{{\mathrm Z}}

\def\bfK{{\mathrm K}}

\def\mo{{\mathfrak o}}

\def\Pic{{\rm Pic}}

\def\PGL{{\rm PGL}}

\def\rk{{\rm rk\,}}

\def\zZ{{\mathcal Z}}
\def\eps{{\epsilon}}

\def\ovl{\overline}

\def\ba{\backslash}

\makeatother
\makeatletter

\author{Joseph Shalika and Yuri Tschinkel}

\title{Height zeta functions of 
equivariant compactifications of the Heisenberg group}

\begin{document} 
 
\date{\today}



\maketitle

\begin{abstract}
We study analytic properties of height zeta functions
of equivariant compactifications of the Heisenberg group.
\end{abstract}

\tableofcontents

\setcounter{section}{0}
\section*{Introduction}
\label{sect:introduction}

Let $\bG=\bG_3$ be the three-dimensional Heisenberg group:
$$
\bG=\{ g=g(x,z,y)=\left(\begin{array}{ccc} 
                                              1 &  x    &  z\\
                                              0 &  1    & y  \\
                                              0 &   0   & 1
                  \end{array}\right) \}. 
$$
Let $X$ be a projective equivariant compactification of $\bG$
(for example $X=\P^3$).  
Thus $X$ is a projective algebraic variety over $\Q$,
equipped with a (left) action of $\bG$  
(and containing $\bG$ as a dense Zariski open subset).
Such varieties can be constructed as
follows: consider a $\Q$-rational 
algebraic representation $\rho\,:\, \bG\ra \PGL_{n+1}$
and take $X\subset \P^n$ to be the Zariski closure 
of an orbit (with trivial stabilizer). This closure need not be smooth
(or even normal). Applying $\bG$-equivariant resolution of singularities
and passing to a desingularization, we may assume that 
$X$ is smooth and that the boundary $D=X\setminus \bG$ 
consists of geometrically irreducible components 
$D=\cup_{\al\in \cA}D_{\al}$, intersecting transversally. 
In this paper, we will always assume that $X$ is a 
bi-equivariant compactification, that is, $X$ carries a left and right
$\bG$-action, extending the left and right 
action of $\bG$ on itself. Equivalently, $X$ is an equivariant
compactification of the homogeneous space $\bG\times \bG/\bG$.

\

Let $L$ be a very ample line bundle on $X$. 
It defines an embedding of $X$ into some projective space $\P^n$.
Let $\cL=(L,\|\cdot \|_{\A})$ be a (smooth adelic) metrization of $L$ and 
$$
H_{\cL}\,:\, X(\Q)\ra \R_{>0}
$$
the associated (exponential) height. 
Concretely, fix a basis $\{ f_j\}_{j=0,...,n}$ 
in the vector space of global sections of $L$ and put
$$
H_{\cL}(x):=\prod_p \max_j(|f_j(x)|_p)\cdot (\sum_{j=0}^n f_j(x)^2)^{1/2}.
$$ 
We are interested in the asymptotics of 
$$
N(B)=N(\cL,B):=\{ \gamma\in \bG(\Q)\,|\, H_{\cL}(\gamma)\le B \}  
$$
as $B\ra \infty$. 

\

The main result of this paper is the determination of 
the asymptotic behavior of $N(B)$ for arbitrary bi-equivariant 
compactifications $X$ of $\bG$ and 
arbitrary projective embeddings. 

\

To describe this asymptotic behavior it is necessary to introduce 
some geometric notions. Denote by $\Pic(X)$ the Picard group of 
$X$. For smooth equivariant compactifications
of unipotent groups, $\Pic(X)$ is freely generated by the classes of
$D_{\al}$ (with $\al\in \cA$). We will use these classes as a basis.
In this basis, the (closed)
cone of effective divisors  $\La_{\rm eff}(X)\subset \Pic(X)_{\R}$
consists of classes 
$$
[L]=(l_{\al})=\sum_{\al\in \cA}l_{\al}[D_{\al}]\in \Pic(X)_{\R},
$$ 
with $l_{\al}\ge 0$ for all $\al$. 
Let $\cL=(L,\|\cdot \|_{\A})$ be a metrized line bundle on $X$ such that 
its class $[L]$
is contained in the interior of $\La_{\rm eff}(X)$.
Conjecturally, at least for varieties
with sufficiently positive anticanonical class,  
asymptotics of rational points of bounded height are related to 
the location of (the class of) $L$ in $\Pic(X)$
with respect to the anticanonical class
$[-K_X]=\kappa=(\kappa_{\al})$ and the cone
$\Lambda_{\rm eff}(X)$
(see \cite{FMT}, \cite{peyre} and \cite{BT}).
In the special case of $\bG$-compactifications $X$ as above and  
$[L]=(l_{\al})$, define:
\begin{itemize}
\item 
$a(L):=\inf\{ a\,\, | \,\, a[L]+[K_X]\in \La_{\rm eff}(X)\}= 
\max_{\al}(\kappa_{\al}/l_{\al});$
\item $b(L):=\#\{ \al\,|\, \kappa_{\al}=a(L)l_{\al}\}$;
\item $\cC(L):=\{ \al \,|\, \kappa_{\al}\neq a(L)l_{\al}\}$;
\item $c(L):=\prod_{\al \notin \cC(L)} l_{\al}^{-1}$.
\end{itemize}
Let 
$$
\zZ(s,\cL):=\sum_{\gamma\in \bG(\Q)} H_{\cL}(\gamma)^{-s},
$$
be the {\em height zeta function} 
(the series converges a priori to a holomorphic function for 
ample $\cL$ and $\Re(s)\gg 0$). 
The Tauberian theorems relate the asymptotics of $N(\cL,B)$  
to analytic properties of $\zZ(s,\cL)$.

\

\begin{Theo}
\label{thm:main}
Let $X$ be a smooth projective bi-equivariant 
compactification of the Heisenberg group $\bG$ 
and $\cL=(L,\|\cdot \|_{\A})$ a line bundle (equipped with a smooth
adelic metrization) such that
its class $[L]\in\Pic(X)$ is contained in the interior of the 
cone of effective divisors $\La_{\rm eff}(X)$.
Then 
$$
\zZ(s,\cL)=\frac{c(L)\tau(\cL)}{(s-a(L))^{b(L)}} +
\frac{h(s)}{(s-a(L))^{b(L)-1}},
$$
where $h(s)$ is a holomorphic function 
(for $\Re(s)>a(L)-\eps$, some $\eps>0$) and $\tau(\cL)$ 
is a positive real number. 
Consequently,
$$
N(\cL,B)\sim\frac{c(L)\tau(\cL)}{a(L)(b(L)-1)!} 
B^{a(L)} \log(B)^{b(L)-1}
$$
as $B\ra \infty$.   
\end{Theo} 

\

\begin{rem}
The  constant $\tau(-\cK_X)$ is the Tamagawa number
associated to the metrization of the anticanonical line bundle
(see \cite{peyre}). For arbitrary polarizations
$\tau(\cL)$ has been defined in \cite{BT}. 
\end{rem}

\

The paper is structured as follows:
in Section~\ref{sect:geom} we
describe the relevant geometric invariants of 
equivariant compactifications of unipotent groups. 
In Section~\ref{sect:metrics} 
we introduce the height pairing
$$
H=\prod_p H_p \cdot H_{\infty}\,:\, \Pic(X)_{\C}\times \bG(\A)\ra \C,
$$
between the complexified Picard group and the adelic points of $\bG$,
generalizing the usual height, and 
the height zeta function
\begin{equation}
\label{eqn:zz}
\zZ({\bf s},g):=\sum_{\gamma\in \bG(\Q)}H({\bf s},\gamma g)^{-1}.
\end{equation}
By the projectivity of $X$, the series 
converges to a function which is continuous and bounded in $g$ and 
holomorphic in ${\bf s}$, for $\Re({\bf s})$ 
contained in some (shifted) cone $\Lambda\subset \Pic(X)_{\R}$.
Our goal is to obtain a meromorphic
continuation of $\zZ({\bf s},g)$ to the tube domain $\sT$
over an open neighborhood of $[-K_X]=\kappa\in \Pic(X)_{\R}$
and to identify the poles.

\

The bi-equivariance of $X$ implies that $H$ is 
invariant under the action {\em on both sides} 
of a compact open subgroup $\bfK$ of the
finite adeles $\bG(\A_{\rm fin})$. Moreover, $H_{\infty}$ is smooth. 
We observe that
$$
\zZ\in \sL^2(\bG(\Q)\backslash \bG(\A))^{\bfK}.
$$
Next, we have, for $\Re({\bf s})$ contained in some
shifted cone in $\Pic(X)_{\R}$, 
an identity in $\sL^2(\bG(\Q)\ba \bG(\A))$ (Fourier expansion):
\begin{equation}
\label{eqn:formalz}
\zZ({\bf s},g)=\sum_{\varrho} \zZ_{\varrho}({\bf s},g),
\end{equation}
where the sum is over all irreducible unitary representations
$(\varrho,\cH_{\varrho})$ of $\bG(\A)$ 
occuring in the right regular representation of $\bG(\A)$ in  
$\sL^2(\bG(\Q)\ba \bG(\A))$ and having $\bfK$-fixed vectors. 
We recall the relevant results
from representation theory in Section~\ref{sect:reps}.

\

We will establish the above identity as an identity of continuous functions
by analyzing the individual terms on the right.
Thus we need to use the (well-known) theory of  
irreducible unitary  representations of the Heisenberg group. 
We will see that for $L=-K_X$ the pole of highest order
of the height zeta function is supplied by the trivial representation.
This need not be the case for other line bundles. Depending on the 
geometry of $X$, it can happen that infinitely many 
non-trivial representations  
contribute to the leading pole of $\zZ(s,\cL)$. 
In such cases the coefficient 
at the pole of highest order is an infinite (convergent) sum of Euler products. 

\

To analyze the contributions in (\ref{eqn:formalz}) 
from the various representations, we 
need to compute local {\em height integrals}. For example,
for the trivial representation, we need to compute the integral  
$$
\int_{\bG(\Q_p)}
H_p({\bf s},g_p)^{-1}dg_p
$$
for almost all $p$ (see Section~\ref{sect:nonar-good}). 
This has been done in \cite{CLT} for
equivariant compactifications of additive groups $\ga^n$; 
the same approach applies here. We regard the height integrals
as geometric versions of Igusa's integrals. They are closely
related to ``motivic'' integrals of Batyrev, Kontsevich, Denef and 
Loeser (see \cite{igusa}, \cite{denef-loeser} and \cite{looijenga}).  

\

\no
The above integral is in fact equal to:
\begin{equation}
\label{eqn:hi}
p^{-\dim(X)}\left( \sum_{A\subseteq \cA} \# D^0_{A}(\bF_p)  
\prod_{\al\in A} \frac{p-1}{p^{s_{\al}-\kappa_{\al}+1}-1} \right),
\end{equation}
where 
$$
D_{\emptyset} := \bG,\,\,\, \, 
D_A:= \cap_{\al\in A}D_{\al},\,\,\, \, 
D_A^0:=D_A\setminus \cup_{A'\supsetneq A}D_{A'}, 
$$
and $\bF_p$ is the finite field $\Z/p\Z$. The resulting Euler product
has a pole of order $\rk\, \Pic(X)$ at ${\bf s}=\kappa$ and 
also the expected leading coefficient at this pole.

\

The bi-$\bfK$-invariance of the height 
insures us that the trivial representation
is ``isolated'' (c.f. especially Proposition~\ref{prop:est-ch}).
Using ``motivic'' integration as above,
we prove that {\em each} of the terms 
on the right side in  (\ref{eqn:formalz})
admits a meromorphic continuation.
We will identify the poles 
of $\zZ_{\varrho}$ for non-trivial representations: 
for ${\bf s}\in \sT$ 
they are contained in the real hyperplanes 
$s_{\alpha}=\kappa_{\alpha}$ and the order of the pole 
at ${\bf s}=\kappa$ is strictly smaller 
than $\rk\, \Pic(X)$. 
Finally, it will suffice to prove the
convergence of the series (\ref{eqn:formalz}), 
for ${\bf s}$ in the appropriate domain. 
This is done in Section~\ref{sect:nonar-good}.

\

This paper is part of a program initiated in \cite{FMT} to relate
asymptotics of rational points of bounded height to geometric invariants.
It continues the work of Chambert-Loir and the second author 
on compactifications of additive groups \cite{CLT}.    
Many statements are direct generalizations from that paper. 
In this paper we explore the interplay between 
the theory of infinite-dimensional representations 
of adelic groups and the theory of height zeta functions
of algebraic varieties. 
The main theorem holds for bi-equivariant compactifications
of arbitrary unipotent groups. We decided to
explain in detail, in a somewhat expository
fashion,  our approach in the simplest possible case
of the Heisenberg group over $\Q$ and 
to postpone the treatment of the general case to a subsequent 
publication.
We have also included the example of $\P^3$ in which most
of the technicalities are absent.

\

{\bf Acknowledgements.}
The second author was partially supported by the NSA, NSF and the 
Clay Foundation.

\

\section{Geometry}
\label{sect:geom}

\begin{nota}
\label{nota:geom}
Let $X$ be a smooth projective algebraic variety. We denote by 
$\Pic(X)$ the Picard group, by $\La_{\rm eff}(X)$ the (closed)
cone of effective divisors and  by $K_X$ the canonical class of $X$. 
If $X$ admits an action by a group $\bG$, we write $\Pic^{\bG}(X)$ for
the group of (classes of) $\bG$-linearized line bundles on $X$. 
\end{nota}

\begin{defn}
\label{defn:ab}
Let $X$ be a smooth projective algebraic variety.
Assume that $\La_{\rm eff}(X)$ is a 
finitely generated polyhedral cone. Let
$L$ be a line bundle such that its class $[L]$ is 
contained in the interior of $\La_{\rm eff}(X)$. 
Define
$$
a(L)=\inf\{ a\,|\, a[L]+[K_X]\in \La_{\rm eff}(X)\}
$$
and $b(L)$ as the codimension of the face of $\La_{\rm eff}(X)$
containing $a(L)[L]+[K_X]$. 
\end{defn}

\begin{nota}
Let $\bG$ be a linear algebraic group over a number field $F$. 
An algebraic variety $X$ (over $F$) will be called a {\em good} 
compactification of $\bG$ if:
\begin{itemize}
\item $X$ is smooth and projective;
\item $X$ contains $\bG$ as a dense Zariski open subset and 
the action of $\bG$ on itself (by left translations) extends to $X$;
\item the boundary $X\setminus \bG$ is a union of 
smooth geometrically irreducible divisors intersecting transversally
(a divisor with strict normal crossings).
\end{itemize}
\end{nota}

\begin{rem}
\label{rem:smooth}
Equivariant resolution of singularities (over a field of characteristic zero)
implies that for {\em any} equivariant compactification 
$X$ there exists an equivariant desingularization 
(a composition of equivariant blowups) 
$\rho\,:\, \tilde{X}\ra X$ such that 
$\tilde{X}$ is a good compactification. 
By the functoriality of heights, 
the counting problem for a metrized line bundle
$\cL$ on $X$ can then be transferred to a counting problem 
for $\rho^*(\cL)$ on $\tilde{X}$. 
Thus it suffices to prove Theorem~\ref{thm:main}
for good compactifications (the answer, of course, does not depend 
on the chosen desingularization). 
\end{rem}

\begin{prop} 
\label{prop:geometry}
Let $X$ be a good compactification of a 
unipotent algebraic group $\bG$. 
Let $D:=X\setminus \bG$ be the boundary and  
$\{D_{\al}\}_{\al\in \cA}$ the set of 
its irreducible components. Then: 
\begin{itemize}
\item  $\Pic^{\bG}(X)_{\Q}=\Pic(X)_{\Q}$;
\item  $\Pic(X)$ is freely generated by the classes $[D_{\al}]$;
\item  $\La_{\rm eff}(X)=\oplus_{\al} \R_{\ge 0} [D_{\al}]$;
\item  $[-K_X]=\sum_{\al}\kappa_{\al}[D_{\al}]$
with $\kappa_{\al}\ge 2 $ for all $\al\in \cA$. 
\end{itemize}
\end{prop}

\begin{proof}
Analogous to the proofs in Section 2 of \cite{HT}.
In particular, it suffices to assume that $X$ carries only 
a one-sided action of $\bG$.     
\end{proof}

\begin{nota}
\label{nota:s}
Introduce coordinates on $\Pic(X)$ using 
the basis $\{ D_{\al}\}_{\al\in \cA}$: 
a vector ${\bf s}=(s_{\al})$
corresponds to $\sum_{\al} s_{\al} D_{\al}$. 
\end{nota}

\begin{coro}
\label{coro:e}
The divisor of every  non-constant function
$f\in F[\bG]$ can be written as 
$$
\dv(f)=E(f)-\sum_{\al} d_{\al}(f)D_{\al},
$$
where $E(f)$ is the unique irreducible component of 
$\{f=0\}$ in $\bG$ and $d_{\al}(f)\ge 0$ for all $\al$. 
Moreover, there is at least one $\al\in  \cA$ such that $d_{\al}(f)>0$. 
\end{coro}

\

\section{Height zeta function}
\label{sect:metrics}

\begin{nota}
For a number field $F$, we denote by $\Val(F)$
the set of all places of $F$, by $S_{\infty}$ the 
set of archimedean and by $S_{\rm fin}$ 
the set of non-archimedean places. For any finite set  $S$ of places
containing $S_{\infty}$, we denote by $\mo_S$ the ring of $S$-integers.
We denote by $\A$ (resp. $\A_{\rm fin}$) the ring of adeles
(resp. finite adeles). 
\end{nota}

\begin{defn}
\label{defn:metri}
Let $X$ be a smooth projective algebraic variety over
a number field $F$. 
A smooth adelic metrization of a line bundle $L$ on $X$
is a family  $\|\cdot \|_{\A} $ of $v$-adic norms 
$\|\cdot \|_v$ on $L\otimes_F{F_v}$,
for all $v\in \Val(F)$, such that:
\begin{itemize}
\item  for $v\in S_{\infty}$ the norm $\|\cdot \|_v$ is  
$\sC^{\infty}$;
\item  for $v\in S_{\rm fin}$, the norm 
of every local section of $L$ is locally constant in 
the $v$-adic topology;
\item  there exist a finite set $S\subset \Val(F)$,
a flat projective scheme (an integral model) 
$\mathcal X$ over $\Spec(\mathfrak o_S)$
with generic fiber $X$ together with a line bundle $\cL$ on 
${\mathcal X}$, such that for all $v\notin S$, the 
$v$-adic metric is given by the integral model.
\end{itemize}
\end{defn}

\begin{prop}
\label{prop:height}
Let $\bG$ be a unipotent
algebraic group defined over a number field $F$ 
and $X$ a good  bi-equivariant compactification of $\bG$.  
Then there exists a height pairing
$$
H=\prod_{v\in \Val(F)} H_v\, :\, \Pic(X)_{\C}\times \bG(\A)\ra \C
$$
such that:
\begin{itemize}
\item 
for all $[L]\in \Pic(X)$, the restriction of 
$H$ to $[L]\times \bG(F)$ 
is a height corresponding to some smooth adelic
metrization of $L$;   
\item 
the pairing is exponential in the $\Pic(X)$ component:
$$
H_v({\bf s}+{\bf s}',g)=H_v({\bf s},g)H_v({\bf s}',g),
$$
for all ${\bf s},{\bf s}'\in \Pic(X)_{\C}$, 
all $g\in \bG(\A)$ and all $v\in \Val(F)$; 
\item there exists a compact open subgroup (depending on $H$) 
$$
\bfK=\bfK(H)=\prod_v \bfK_v\subset \bG(\A_{\rm fin})
$$ 
such that, for all $v\in S_{\rm fin}$,  one 
has $H_v({\bf s},kgk')=H_v({\bf s} ,g)$ 
for all ${\bf s}\in \Pic(X)_{\C}$, $k,k'\in \bfK_v$ and $g\in \bG(F_v)$.  
\end{itemize}
\end{prop}

\begin{proof}
For $\bG=\ga^n$ the Proposition is proved in \cite{CLT}, Lemma 3.2. 
The same proof applies to any unipotent group.
\end{proof}

\begin{nota}
\label{nota:cone}
For $\delta \in \R$, we denote by 
$\sT_{\delta}\subset \Pic(X)_{\C}$ the tube domain 
$\Re({s}_{\al})-\kappa_{\al} > \delta$ (for all $\al \in \cA$). 
\end{nota}

\begin{defn}
\label{defn:zeta}
The height zeta function on $\Pic(X)_{\C}\times \bG(\A)$ is defined as
$$
\zZ({\bf s},g)=\sum_{\gamma\in \bG(F)} H({\bf s},\gamma g)^{-1}.
$$ 
\end{defn}

\begin{prop}
\label{prop:abs-conv}
There exists a  $\delta >0$ 
such that, for all ${\bf s}\in \sT_{\delta}$ and all $g\in \bG(\A)$, 
the series defining the height zeta function $\zZ({\bf s},g)$
converges normally
(for $g$ and ${\bf s}$ contained in compacts in
$\bG(\A)$, resp. $\sT_{\delta}$)
to a function which is holomorphic in ${\bf s}$ and continuous in $g$.
\end{prop}

\begin{proof}
The proof is essentially analogous 
to the proof of Proposition 4.5 in \cite{CLT}
(and follows from the projectivity of $X$).
\end{proof}

\begin{coro}
\label{coro:formal}
For ${\bf s}\in \sT_{\delta}$, one has  
an identity in $\sL^2(\bG(F)\ba \bG(\A))$, as above:
\begin{equation}
\label{eqn:formal}
\zZ({\bf s},g)=\sum_{\varrho} \zZ_{\varrho}({\bf s},g).
\end{equation}
The sum is over all irreducible unitary representations 
$\varrho $ of $\bG(\A)$ occuring  $\sL^2(\bG(F)\ba \bG(\A))$ 
and having a $\bfK$-fixed vector (cf. Proposition~\ref{prop:formm}). 
\end{coro}

\section{Representations}
\label{sect:reps}

\subsection{}
\label{sect:30}

From  now on, for the sake of simplicity,  we suppose $F=\Q$. 
Denote by $\bZ=\ga$ the one-dimensional center  
and by $\bG^{\rm ab}=\bG/\bZ=\ga^{2}$ the
abelianization of $\bG$. 
Let $\bU\subset \bG$ be the subgroup
$$
\bU:=\{ u\in \bG\,|\, u=(0,z,y)\}
$$
and  
$$
\bW:= \{ w\in \bG\,|\, w=(x,0,0)\}.
$$
We have $\bG=\bW\cdot \bU=\bU\cdot \bW$. We may assume that 
the compact open subgroup
$$
\bfK=\prod_p \bfK_p \subset \bG(\A_{\rm fin})
$$ 
of Proposition~\ref{prop:height}
is given by
\begin{equation}
\label{eqn:k}
\bfK=\prod_{p\notin S_{H}} \bG(\Z_p) \cdot 
\prod_{p\in S_{H}} \bG(p^{n_p}\Z_p),
\end{equation}
where $S_H$ is a finite set of primes and the $n_p$ are positive 
integers. 
We denote by $\bfK^{\rm ab},\bfK_{\bZ}$ etc. 
the corresponding compact subgroups of the (finite)
adeles of $\bG^{\rm ab},\bZ, \bU, \bW$, respectively, and 
put 
$$
\sn(\bfK)=\prod_{p\in S_H} p^{n_p}.
$$
We denote by $dg=\prod_p dg_p \cdot dg_{\infty}$ the Haar measure
on $\bG(\A)$, where we have set $dg_p=dx_pdy_pdz_p$ with the normalization 
$\int_{\Z_p} dx_p = 1$ etc. (similarly at the real place). 
We write $du_p=dz_pdy_p$ (resp. $du_{\infty}$, $du$) for the Haar measure 
on $\bU(\Q_p)$ (resp. $\bU(\R)$, $\bU(\A)$).  
We let $dk_p$ be the Haar measure on $\bfK_p$ obtained by 
restriction of $dg_p$ to $\bfK_p$. Further, our normalization of
measures implies that $\int_{\bfK_p} dk_p =1$. 
As usual, a choice of a measure on the local (or global) points of 
$\bG$ and of a subgroup $\bH\subset \bG$
determines a unique measure on the local 
(resp. global) points of the homogeneous space
$\bG/\bH$.

\begin{lemm}
\label{lemm:K}
One has:
\begin{itemize}
\item $\bG(\Z_p)
=(\bG(\Z_p)\cap \bU(\Q_p))\cdot (\bG(\Z_p)\cap \bW(\Q_p))$;
\item $\bU(\Q_p)\cdot \bW(\Z_p)$ is a subgroup of $\bG(\Q_p)$;
\item $\bG(\A)=\bG(\Q)\cdot \bG(\R)\cdot \bfK;$
\item there exists a subgroup $\Gamma \subset \bG(\Z)$ (of finite index) 
such that 
$$
\bG(\Q)\backslash \bG(\A)/\bfK 
= \Gamma \backslash \bG(\R);
$$
\item the quotient $\Gamma\ba \bG(\R)$ is compact.
\end{itemize}
\end{lemm}
These statements are well-known and easily verified.

\

We now recall the well-known representation theory of the 
Heisenberg group in an adele setting (\cite{howe}).
Denote by  $\varrho$ the right regular
representation of $\bG(\A)$  on the Hilbert space
$$
\cH:=\sL^2(\bG(\Q)\ba\bG(\A)).
$$
 
Consider the action of the compact group 
$\bZ(\A)/\bZ(\Q)$ on $\cH$ (recall that $\bZ=\ga$). 
By the Peter-Weyl theorem, we obtain a decomposition
$$
\cH=\oplus \cH_{\psi}
$$
and corresponding representations 
$(\varrho_{\psi},\cH_{\psi})$ of $\bG(\A)$. Here
$$
\cH_{\psi}:=\{ \varphi\in \cH \, |\, 
\varrho(z)(\varphi)(g)=\psi(z)\varphi(g) \}
$$
and $\psi$ runs over the set of (unitary)
characters of $\bZ(\A)$ which are trivial on $\bZ(\Q)$. 
For non-trivial $\psi$, the corresponding representation
$(\varrho_{\psi},\cH_{\psi})$ of $\bG(\A)$ 
is non-trivial, irreducible and unitary. 
On the other hand, when $\psi$ is the trivial character, 
the corresponding representation
$\varrho_{0}$ decomposes further as a 
direct sum of one-dimensional representations $\varrho_{\eta}$:
$$
\cH_0=\oplus_{\eta} \cH_{\eta}.
$$
Here $\eta$ runs once over all (unitary) characters of the group 
$\bG^{\rm ab}(\Q)\ba \bG^{\rm ab}(\A)$. 
It is convenient to consider 
$\eta$ as a function on $\bG(\A)$, 
trivial on the $\bZ(\A)$-cosets.
Precisely, let $\psi_1=\prod_p \psi_{1,p} \cdot \psi_{1,\infty}$ 
be the Tate-character 
(which has exponent zero at each finite prime, see \cite{tate67b} 
and \cite{weil}).
For $\mathbf a=(a_1,a_2)\in\A\oplus \A$, consider
the corresponding linear form on 
$$
\bG^{\rm ab}(\A)=\A\oplus \A
$$ 
given by 
$$
g(x,z,y)\mapsto a_1x+a_2y
$$ 
and denote by $\eta=\eta_{\mathbf a}$ (${\bf a}=(a_1,a_2)$) the corresponding
adelic character 
$$
\eta\,:\, g(x,z,y))\mapsto \psi_1(a_1x+a_2y)
$$ 
of $\bG(\A)$.
For $a\in \A$, we will denote by $\psi_a$ the adelic character 
of $\bZ(\A)$ given by 
$$
z\mapsto \psi_1(az).
$$ 

\

As in Section~\ref{sect:metrics}, 
the starting point of our 
analysis of the height zeta function is the spectral decomposition of $\cH$. 
A more detailed version of Corollary~\ref{coro:formal} is
the following Proposition.

\begin{prop}
\label{prop:formm}
There exists a $\delta>0$ such that,
for all ${\bf s}\in \sT_{\delta}$, one has 
an identity of $\sL^2$-functions
\begin{equation}
\label{eqn:id}
\zZ({\bf s},g)=\zZ_0({\bf s},g)+\zZ_1({\bf s},g)+\zZ_2({\bf s},g),
\end{equation}
where
\begin{equation}
\label{eqn:zo}
\zZ_0({\bf s},id)  =  \int_{\bG(\A)}H({\bf  s},g)^{-1}dg,
\end{equation}
\begin{equation}
\label{eqn:z1}
\zZ_1({\bf s},g) = \sum_{\eta} \eta(g)\cdot  \zZ({\bf s},\eta),
\end{equation}
and 
\begin{equation}
\label{eqn:z2}
\zZ_2({\bf s},g) =  \sum_{\psi} 
\sum_{\omega^{\psi}} \omega^{\psi}(g)\cdot
\zZ({\bf s},\omega^{\psi}). 
\end{equation}
Here we have set
$$
\zZ({\bf s},\eta):=\langle \zZ({\bf s},\cdot),\eta\rangle=
\int_{\bG(\A)}H({\bf s},g)^{-1}\overline{\eta}(g)dg,
$$
$$
\zZ({\bf s},\omega^{\psi}):=\langle\zZ(s,g),\omega^{\psi}\rangle=
\int_{\bG(\Q)\backslash \bG(\A)}
\zZ({\bf s},g)\overline{\omega}^{\psi}(g)dg=
$$
$$
 =\int_{\bG(\A)}
H({\bf s},g)^{-1}\overline{\omega}^{\psi}(g)dg,
$$
$\eta$ 
ranges over all non-trivial characters of 
$$
\bG^{\rm ab}(\Q)\cdot {\bfK}^{\rm ab}\ba \bG^{\rm ab}(\A),
$$
$\psi$ ranges over all non-trivial characters of 
$$
\bZ(\Q)\cdot {\bfK}_{\bZ}\ba \bZ(\A),
$$
and $\omega^{\psi}$ ranges over  
a fixed orthonormal basis of $\cH_{\psi}^{\bfK}$
(for each $\psi$).

In  particular, for $\eta=\eta_{\bf a}$ and $\psi=\psi_{a}$
occuring in this decomposition, we have 
$$
a_1,a_2,a\in \frac{1}{\sn(\bfK)}\Z.
$$ 
\end{prop}

\begin{proof}
We use the (right) $\bfK$-invariance of the height for the 
last statement (for $\eta$). For $\psi$ see also Lemma~\ref{lemm:j}
as well as Proposition~\ref{prop:abs-conv}. 
\end{proof}

\begin{rem}
\label{rem:zzz}
The desired meromorphic properties of $\zZ_0$ and $\zZ_1$ 
have, in fact, already been established in \cite{CLT}.
The height integrals are computed
as in the abelian case and 
the convergence of the series $\zZ_1$ is 
proved in the same way as in \cite{CLT}.    
In particular, (\ref{eqn:z1}) is an identity of continuous functions.
The novelty here is the treatment of $\zZ_2$. 
\end{rem}

We now proceed to describe the various standard models
of infinite-dimensional representations
of the Heisenberg group. 

\

\subsection{}
\label{sect:31}

\no
{\bf Locally:}
Let $\psi=\psi_p$ (resp. $\psi=\psi_{\infty}$)
be a {\em local} non-trivial character of $\Q_p$ (resp. $\R$). 
Extend $\psi$ to $\bU(\Q_p)$ by setting
$$
\psi((0,z,y))=\psi(z).
$$
The one-dimensional representation of $\bU(\Q_p)$ thus obtained
induces a representation $\pi_{\psi}=\pi_{\psi,p}$ of $\bG(\Q_p)$.
The representation $\pi_{\psi}$ 
acts on the Hilbert space  
of measurable functions 
$$
\phi\,:\, \bG(\Q_p)\ra \C
$$ 
which satisfy the conditions:
\begin{itemize}
\item 
$\phi(ug)=\psi(u)\phi(g)$ 
for all $u\in \bU(\Q_p)$ and $g\in \bG(\Q_p)$;
\item  
$\|\phi\|^2:=\int_{\bU(\Q_p)\ba\bG(\Q_p)}|\phi(g)|^2 dg<\infty$.
\end{itemize}
The action is given by
$$
\pi_{\psi}(g')\phi(g)=\phi(gg'), \,\,\, g'\in \bG(\Q_p).
$$

On the other hand, we have 
a representation
$\pi_{\psi}'=\pi_{\psi,p}'$ (the oscillator representation) on 
$$
\sL^2(\bW(\Q_p))=\sL^2(\Q_p),
$$ 
where the action of $\G(\Q_p)$ on a function
$\varphi\in \sL^2(\Q_p)$ is given by

\begin{equation}
\label{eqn:10}
\pi_{\psi}' (g(x',0,0)) \varphi(x) = \varphi(x+x')
\end{equation}
$$
\pi_{\psi}' (g(0,0,y))  \varphi(x)  = \psi(y\cdot x)\varphi(x)     
$$
$$
\pi_{\psi}' (g(0,z,0)) \varphi(x)  = 
\psi(z)\varphi(x).
$$

It is easy to see that the representations $\pi_{\psi}$ and $ \pi_{\psi}' $
are unitarily equivalent. We will identify the unitary representations
$\pi_{\psi}$ and $\pi_{\psi}'$ in  what follows.

\

\no
{\bf Globally:}
In the adelic situation, to each non-trivial character 
$\psi$ of $\bZ(\A)$
we can associate a  representation  
$\pi_{\psi}$  of $\bG(\A)$, where
$\pi_{\psi}= \otimes_p \pi_{\psi,p} \otimes \pi_{\psi,{\infty}}$ 
and  the action on $\sL^2(\bU(\A)\ba \bG(\A))=\sL^2(\A)$, 
is given by the formulas (\ref{eqn:10}) 
(with $\psi_p$ replaced by $\psi$). 
The representations  
$\pi_{\psi}$ and $\varrho_{\psi}$
are equivalent irreducible unitary  representations of $\bG(\A)$.
We will recall the explicit 
intertwining map between $\pi_{\psi}$ and $\varrho_{\psi}$
(c.f. Lemma~\ref{lemm:j}).

\

We also recall that the space $\cS(\A)\subset \sL^2(\A) $
of Schwartz-Bruhat functions  
coincides with the space
of smooth vectors of $\pi_{\psi}$ (for the real place, 
see the Appendix in \cite{cg}) and note that 
$\sL^2(\Q_p)^{\bfK_p}=\cS(\Q_p)^{\bfK_p}$.

\

For a character $\psi(z)=\psi_{\infty}(z)=e^{2\pi i az}$
(with $a\neq 0$)  
consider the following operators on the 
subspace of Schwartz functions
$\cS(\R)\subset \sL^2(\R)$:
$$
\begin{array}{rcl}
{\mathfrak d}^+_\psi \varphi(x)  &  = & 
\frac{d}{d x} \varphi(x)\\
{\mathfrak d}^{-}_\psi\varphi(x) & = & 
2\pi i ax\varphi(x)\\
\Delta_{\psi} & = & ({\mathfrak d}^{+}_{\psi})^2+({\mathfrak d}^{-}_{\psi})^2.
\end{array}
$$
We have
$$
\Delta_{\psi}\varphi(x)=\varphi''(x)-(2\pi a x)^2\varphi(x)
$$
(harmonic oscillator). The eigenvalues of $\Delta_{\psi}$ are
given by
$$
\la_n^{\psi} = -2\pi(2n+1)|a|
$$ 
(with $n=0,1,2,...$). They have multiplicity one. 
Denote by $h_n^{\psi}(x)$ the $n$-th
Hermite polynomial: 
$$
\begin{array}{rcl}
h_0^{\psi}(x) & = & 1\\
h_1^{\psi}(x) & = & 4\pi|a|x\\
h_2^{\psi}(x) & = & -4\pi|a|(1-4\pi|a|x^2)
\end{array}
$$
and, in general,
$$
\frac{d^n}{dx^n} e^{-2\pi|a|x^2}=(-1)^nh_n^{\psi}(x)e^{-2\pi|a|x^2}.
$$
The (essentially unique) 
eigenfunction $\varphi^{\psi}_n$ corresponding to $\la_n^{\psi}$
is given by 
$$
\varphi^{\psi}_n:=c_ne^{-\pi|a|x^2}h_n^{\psi}(x).
$$
Here we choose the constants $c_n$, so that the $\sL^2$-norm of
$\varphi^{\psi}_n$ is $1$.

\begin{lemm}
\label{lemm:basis}
The set $\cB_{\infty}(\pi_{\psi}'):=\{\varphi^{\psi}_n\}$ 
is a complete orthonormal basis of $\sL^2(\R)$. 
\end{lemm}

\begin{proof}
For details see, for example,  \cite{b}, Chapter 13,  or \cite{ch}.
\end{proof}

\

\subsection{}
\label{sect:32}

\begin{nota}
\label{nota:eta}
For $\eta=\eta_{\bf a}$ with 
${\bf a}=(a_1,a_2)$ and $a_1,a_2\in \frac{1}{\sn(\bfK)}\Z$, 
denote by $S_\eta$ the set of primes $p$ dividing either 
$\sn(\bfK)a_1$ or $\sn(\bfK)a_2$. 
Similarly, for $\psi=\psi_a$ with $a\in  \frac{1}{\sn(\bfK)} \Z$, 
denote by $S_{\psi}$ the set of primes dividing $\sn(\bfK)a$. 
\end{nota}

\begin{lemm}
\label{lemm:mult}
Let $\psi=\psi_a$ be a non-trivial character of $\bZ(\Q)\ba \bZ(\A)$ and 
$\varrho_{\psi}=\otimes_p \varrho_{\psi,p}\otimes \varrho_{\psi,\infty}$ 
the corresponding infinite-dimensional automorphic representation. 
Suppose $\varrho_{\psi}$ contains a $\bfK$-fixed vector
(for $\bfK$ as in (\ref{eqn:k})). Then: 
\begin{itemize} 
\item
$a\in \frac{1}{\sn(\bfK)}\Z$ (for $\sn(\bfK)=\prod_{p\in S_H}p^{n_p}$);
\item $\dim \varrho_{\psi,p}^{\bfK_p}=1$ for $p\notin S_{\psi}$; 
\item $\dim \varrho_{\psi,p}^{\bfK_p}= |\sn(\bfK)^2a|_p^{-1}$ for $p\in S_{\psi}$, 
provided $p^{n_p}\cdot \sn(\bfK) \in \Z_p$.
\end{itemize}
\end{lemm}

\begin{proof}
We need only use the explicit form of the representation
$\pi_{\psi,p}$ given in (\ref{eqn:10}). Suppose first that
$\pi_{\psi,p}$ has a non-zero $\bfK_p$-fixed vector $\varphi$. 
Taking $z\in p^{n_p}\cdot \Z_p$, we get
$$
\psi_{p}(p^{n_p}r) =\psi_{1,p}(ap^{n_p}r)=1
$$
for all $r\in \Z_p$. Since the exponent of $\psi_{1,p}$ is zero, we have
$$
a\cdot p^{n_p} \in \Z_p,
$$
from which the first assertion follows.

Let us assume then that $p^{n_p}\cdot \sn(\bfK)\in \Z_p$. Then the space of $\bfK_p$-fixed
vectors $\varphi $ in $\sL^2(\Q_p)$ is precisely the set of $\varphi$
satisfying
\begin{enumerate}
\item $\varphi(u+p^{n_p}r_1)=\varphi(u)$;
\item $\varphi(u) =\psi(p^{n_p} r_2 u)\varphi(u)$
\end{enumerate}
for all $r_1,r_2\in \Z_p, u\in\Q_p$. The first identity implies that
$\varphi$ is a continuous function and the second 
that ${\rm Supp}(\varphi) \subset a^{-1} p^{-n_p}\cdot \Z_p$. The second 
and the third assertions of the Lemma follow at once.  
\end{proof}

\begin{nota}
\label{nota:v}
Let $V_{\psi,p}$ be the space of the induced representation of $\pi_{\psi,p}$.
Denote by $V_{\psi,p}^{\infty}$ the space of smooth vectors in $V_{\psi,p}$. 
Thus $V_{\psi,p}^{\infty}$ is the set of all $\nu\in V_{\psi,p}$ fixed by 
some open compact subgroup of $\bG(\Q_p)$. Note that $V_{\psi,p}^{\infty}$ 
is stable under the action of $\bG(\Q_p)$. Note also that in the explicit
realization of $\pi_{\psi,p}$ given in (\ref{eqn:10}), $\sL^2(\Q_p)^{\infty}=\cS(\Q_p)$
(see the proof of Lemma~\ref{lemm:mult}).

\end{nota}

\

\noindent
For $\varphi\in \cS(\A)$ define the theta-distribution
$$
\Theta(\varphi):=\sum_{x \in \Q} \varphi(x).
$$
Clearly, $\Theta$ is a 
$\bG(\Q)$-invariant linear functional on $\cS(\A)$. 
This gives a map
$$
\begin{array}{rcl}
j_{\psi}\,:\, \cS(\A)     &  \ra & \sL^2(\bG(\Q)\ba \bG(\A))\\
      j_{\psi}(\varphi)(g)&  =   & \Theta(\pi_{\psi}(g)\varphi).
\end{array}
$$

\begin{lemm}
\label{lemm:j}
The map $j_{\psi}$ extends to an isometry 
$$
j_{\psi}\,:\, \sL^2(\A)\stackrel{\sim}{\longrightarrow}
\cH_{\psi}\subset \sL^2(\bG(\Q)\ba\bG(\A)),
$$ 
intertwining $\pi_{\psi}$ and $\varrho_{\psi}$.
Moreover,  
$$
j_{\psi}\,:\, 
\sL^2(\A)^{\bfK}\stackrel{\sim}{\longrightarrow}  
\cH_{\psi}^{\bfK}.
$$
\end{lemm}

\

Let us recall the definition of a restricted algebraic tensor product:
for all primes $p$, let $V_p$ be a (pre-unitary) 
representation space for $\bG(\Q_p)$. 
Let $(e_p)_p$ be
a family of vectors $e_p\in V_p$, defined for all primes $p$ outside a finite
set $S_0$. Suppose that, for almost all $p$, $e_p$ is fixed 
by $\bfK_p$. 
We will also assume that the norm of $e_p$ is equal to $1$. 
Let $S$ be a finite set of primes containing $S_0$. 
A {\em pure} tensor is a vector,
$\nu=\nu_S\otimes e^S$, where $e^S=\otimes_{p\notin S}e_p$ and $\nu_S$
is a pure tensor in the finite tensor product  
$\otimes_{p\in S} V_p$. 
The restricted algebraic tensor product $V=\otimes_p V_p$ is 
generated by finite linear combinations of pure tensors
(see \cite{JL} for more details). 

\

\begin{exam}
\label{exam:ep}
Consider the representation 
$\pi_{\psi}$ of $\bG(\A)$ on the 
Schwartz-Bruhat space $\cS(\A_{\rm fin})=\otimes_p \cS(\Q_p)$
and the corresponding representation $\pi_{\psi_p}$ of $\bG(\Q_p)$ 
on $\cS(\Q_p)$.  
In this case, for all primes $p\notin S_{\psi}$, 
$e_p$ is unique (up to scalars) and may be taken to be 
the characteristic function of $\Z_p$.
We have 
$j_{\psi}(\cS(\A_{\rm fin})\otimes \cS(\R))=\cH^{\rm smooth}_{\psi}$ 
(by \cite{cg}). 
\end{exam}

\

We now fix an orthonormal basis 
$\cB_{\rm fin}(\pi_{\psi})$
for the space $\cS(\A_{\rm fin})^{\bfK}$ as follows. 
We let $\cB_{\rm fin}(\pi_{\psi})=\otimes_p \cB_p(\pi_{\psi_p})$, 
where, for $p\in S_0=S=S_{\psi}$, $\cB_p(\pi_{\psi_p})$ is any fixed orthonormal 
basis for $\cS(\Q_p)^{\bfK_p}$ and, for $p\notin S$, 
$\cB_p(\pi_{\psi_p})=e_p$. 
Thus any $\varphi \in \cB_{\rm fin}(\pi_{\psi})$ 
has the form
$$
\varphi=\varphi_S\otimes e^S,
$$ 
with $e^S=\otimes_{p\in S} e_p$, as above.
We have then the following Lemma:

\begin{lemm}
\label{lemm:dimension}
The set 
$$
\cB(\varrho_{\psi}):= j_{\psi}
(\cB_{\rm fin}(\pi_{\psi})\otimes\cB_{\infty}(\pi_{\psi}))
$$
is a complete orthonormal basis of $\cH^{\bfK}_{\psi}$. 
The number of elements $\omega\in \cB(\varrho_{\psi})$
(c.f. Lemma~\ref{lemm:mult})
with given eigenvalue $\la_n^{\psi}$ 
is $|\sn(\bfK)^2a|$ 
if $a\in \frac{1}{\sn(\bfK)}\Z$ (and  zero otherwise). 
\end{lemm}

\begin{defn}
\label{def:sph}
Suppose $p\notin S_{\psi}$. The normalized
spherical function $f_p$ on $\bG(\Q_p)$ is defined by 
$$
f_p(g_p):=\langle \pi_{\psi_p} (g_p)e_p,e_p\rangle.
$$
Here $\langle \cdot ,\cdot \rangle $ is the standard
inner product on $\sL^2(\Q_p)$. 
\end{defn}

\begin{lemm}[Factorization]
\label{lemm:fact}
For $\omega \in \cB(\varrho_{\psi})$ and $S=S_{\psi}\cup \{ \infty\}$, 
we have
an identity 
$$
\int_{\bfK^S} \omega (k^Sg)dk^S =\prod_{p\notin S} f_p(g_p) \cdot \omega(g_{S}).
$$
Here $\bfK^S=\prod_{p\notin S_{\psi}} \bfK_p$, 
$g=g^S\cdot g_S$,  with $g^S$ (resp. $g_S$) in $\bG(\A^S)$ (resp. $\bG(\A_S)$).  
\end{lemm}

\begin{proof}
Define a linear form $\mu$ on $V=\cS(\A)$ by setting
$$
\mu(\varphi):= \int_{\bfK^S} j(\varphi)(k^S)dk^S,
$$
(where $\varphi\in \cS(\A)$).
Set 
$$
V^S:=\otimes_{p\notin S} \cS(\Q_p)
$$
and 
$$
V_S:=\otimes_{p\in S_{\psi}} \cS(\Q_p) \otimes \cS(\R),
$$
so that $V=V_S\otimes V^S$. Then from Lemma~\ref{lemm:mult} we have,
for $\varphi^S\in V^S$, 
with $\pi_{\psi}^S =\otimes_{p\notin S}\pi_{\psi,p}$, an 
equality of the form
$$
 \int_{\bfK^S} \pi_{\psi}^S(k^S)\varphi^S dk^S =\nu^S (\varphi^S)\cdot e^S
$$ 
for a unique linear form $\nu^S$ on $V^S$. Note that
$\nu^S(\varphi^S)=\langle \varphi^S,e^S\rangle$, for $\varphi^S\in V^S$. 
Now we have, for $\varphi $ of the form $\varphi=\varphi_S\otimes \varphi^S$, 
with $\varphi_S\in V_S$ and $\varphi^S\in V^S$, 
$$
\mu(\varphi_S\otimes \pi_{\psi}^S(k^S)\varphi^S) =\mu (\varphi_S\otimes \varphi^S),
$$
from which it follows at once that
$$
\mu(\varphi_S \otimes \varphi^S) =\mu_S(\varphi_S)\cdot \nu^S(\varphi^S),
$$
for some linear form $\mu_S$ on $V_S$.
From this we obtain in turn, for $\varphi=\varphi_S\otimes e^S$, the identity
$$
\begin{array}{ccl}
\int_{\bfK^S} j(\varphi) (k^S g) dk^S  & = & \mu(\pi_{\psi}(g)\varphi) \\
 &   & \\
 & = &  \mu_S(\pi_{\psi,S}(g_S) \varphi_S) \cdot \nu^S (\pi_{\psi}^S(g^S) e_S) \\
 &  &  \\
 & = & \mu_S (\pi_{\psi,S}(g_S) \varphi_S) \cdot \prod_{p\notin S} f_p(g_p)
\end{array}
$$
for $g\in \bG(\A)$. 
Here $\pi_{\psi,S} =\otimes_{p\in S}\pi_{\psi,p}$. 
Taking $\omega =j(\varphi)$, with $\varphi = 
\varphi_S\otimes e^S, \,\, \varphi_S\in V_S$ as above, 
we arrive next at the equality
$$
\int_{\bfK^S} \omega(k^Sg) dk^S = \omega'(g_S)\cdot \prod_{p\notin S}f_p(g_p),
$$
for some function $\omega'$ on $\bG(\A_S)$. Finally, if $g=g_S\in \bG(\A_S)$, 
we obtain from the last expression
$$
\omega'(g_S) =\int_{\bfK^S} \omega(k^Sg_S)dk^S = 
\int_{\bfK^S}\omega(g_Sk^S)dk^S =\omega(g_S), 
$$
since, in fact, $\omega$ is $\bfK$-invariant on the right. This completes the proof
of the Lemma. 
\end{proof}

\

\begin{coro}
\label{coro:fact}
Let $\psi=\psi_a$ be as above (with $a\in \Q^\times$) and  
$\varrho_{\psi}$  the associated
irreducible unitary automorphic representation of $\bG(\A)$. 
Suppose that
$\varrho_{\psi}$ has a $\bfK$-fixed vector.   
Then, for $S=S_{\psi}\cup \{ \infty\}$, 
all $\omega\in \cB(\varrho_{\psi})$, all primes 
$p\notin S_{\psi}$ and all (integrable)
functions $H$ on $\bG(\A)$ 
such that 
$$
H_p(k_pg_p)=H_p(g_pk_p)=H_p(g_p),
$$ 
for all $k_p\in \bfK_p$ and $g_p\in \bG(\Q_p)$, 
one has 
\begin{equation}
\label{eqn:11}
\int_{\bG(\A)}H(g)\omega(g)dg=
\prod_{p\notin S} \int_{\bG(\Q_p)}  H_p(g_p) f_p (g_p) dg_p \cdot
\int_{\bG(\A_S)} H(g_S)\omega_S (g_S) dg_S,   
\end{equation}
where $\omega_S$ is the restriction of $\omega$ to $\bG(\A_S)$. 
\end{coro}

\begin{lemm}
\label{lemm:u}
For all $\psi$ and all $p\notin S_{\psi}$ 
one has, for $H_p$ as above, 
$$
\int_{\bG(\Q_p)}  H_p(g_p) f_p (g_p) dg_p = 
\int_{\bU(\Q_p)} H_p(u_p)\psi_p(u_p)du_p.
$$
\end{lemm}

\begin{proof}
Suppose $p\notin S_{\psi}$. Let $\chi_p$ be the characteristic function
of $\bfK_p$. Define a function $\tilde{\psi}_p$ on $\bG(\Q_p)$ by setting
$$
\tilde{\psi}_p(g_p):= \int_{\bU(\Q_p)} \chi_p (u_pg_p)\ovl{\psi}_p(u_p)du_p
$$
(with $g_p\in \bG(\Q_p)$). Clearly, $ \tilde{\psi}_p $ belongs to the space
$V_{\psi,p}$ of the induced representation $\pi_{\psi,p}$; moreover, 
$\tilde{\psi}_p$ is $\bfK_p$-invariant (on the right). 

Next we have, with our normalization of Haar measures, 
$$
\tilde{\psi}_p(g_p) = \psi_p(u_p)
$$
provided $g_p=u_p k_p$, with $u_p\in \bU(\Q_p), k_p\in \bfK_p$, and zero otherwise. 
In particular, 
$$
|\tilde{\psi}_p(g_p)|^2 =\int_{\bU(\Q_p)} \chi_p(u_pg_p)du_p,
$$ 
from which it follows that
$$
\| \tilde{\psi}_p\|^2=\int_{\bU(\Q_p)\ba \bG(\Q_p)} |\tilde{\psi}_p(g_p)|^2 d^*g_p 
= \int_{\bG(\Q_p)} \chi_p(g_p)dg_p =\int_{\bfK_p}dg_p =1.
$$
(Here $d^*g_p$ is normalized so 
that $dg_p=du_p d^*g_p$ as in Section~\ref{sect:30}.)
Next, for $\nu\in V_{\psi,p}^{\infty}$, we have, with $\pi_p=\pi_{\psi,p}$, 
$$
\int_{\bfK_p}\pi_p(k_p)\nu dk_p=\mu(\nu)\tilde{\psi}_p,
$$
for a unique linear form $\mu$ on $V_{\psi,p}^{\infty}$. Note that
$\mu(\nu)=\langle \nu,\tilde{\psi}_p\rangle$. 
We have then, using $\tilde{\psi}_p(e)=1$, 
$$
f_p(g_p) =  \langle \pi_p(g_p)\tilde{\psi}_p,\tilde{\psi}_p\rangle =
\int_{\bfK_p} \tilde{\psi}_p(k_pg_p)dk_p. 
$$
To complete the proof, we note first, from the left $\bfK_p$-invariance
of $H_p$, that 
$$
\int_{\bG(\Q_p)} H_p(g_p)f_p(g_p)dg_p = \int_{\bG(\Q_p)} H_p (g_p) \tilde{\psi}_p (g_p) dg_p.
$$ 
In turn, the last integral is
$$
\begin{array}{ccl}
  & = & \int_{\bU(\Q_p)} \ovl{\psi}_p(u_p) \int_{\bG(\Q_p)} H_p(g_p)\chi_p(u_pg_p) dg_p\\
 &  & \\
& = & \int_{\bU(\Q_p)} \psi_p(u_p) \int_{\bG(\Q_p)}H_p(u_pg_p)\chi_p(g_p)dg_p\\
  & &  \\
& = & \int_{\bU(\Q_p)} H_p(u_p)\psi_p (u_p)du_p,
\end{array}
$$
the last equality from  the right $\bfK_p$-invariance of $H_p$. 

\end{proof}

\section{Euler Products}
\label{sect:nonar-good}

In this section we show that each summand in 
the $\sL^2$-expansion of the height zeta function
in Proposition~\ref{prop:formm} 
is regularized by an explicit Euler product. 
First we record the integrability of local heights:

\begin{lemm}
\label{lemm:bad-red}
For all compacts $\sK\subset \sT_{-1}$
and all primes $p$,
there exists a constant $\sc_p(\sK)$
such that, for all ${\bf s}\in \sK$,
one has:
$$
\int_{\bG(\Q_p)}|H_p({\bf s},g_p)^{-1}|dg_p
\le \sc_p(\sK).
$$
Moreover, for all 
$\partial$ in the universal enveloping algebra 
$ \mathfrak U(\mathfrak g)$ of $\bG$ and all compacts 
$\sK\subset \sT_{-1}$,
there exists a constant $\sc(\sK,\partial)$
such that, for all ${\bf s}\in \sK$,
$$
\int_{\bG(\R)}|\partial H_\infty({\bf s},g_\infty)^{-1}|dg_\infty
\le \sc(\sK,\partial).
$$
\end{lemm}

\begin{proof}
This is Lemma 8.2 and Proposition 8.4 of \cite{CLT}.
\end{proof}

\begin{nota}
\label{nota:SX}
Denote by $S_X$ the set of all primes such that
one of the following holds:
\begin{itemize}
\item $p$ is $2$ or $3$;
\item $\bfK_p\neq \bG(\Z_p)$;
\item over $\Z_p$, the union $\cup_{\al}D_{\al}$ is not
a union of smooth relative divisors with strict normal crossings.
\end{itemize}
\end{nota}

\begin{rem}
\label{rem:inva}
For all $p\notin S_X$, 
the height $H_p$ is invariant with respect to
the right and left $\bG(\Z_p)$-action. 
\end{rem}

\begin{prop}
\label{thm:trivial}
For all primes $p\notin S_X$ and all ${\bf s}\in 
\sT_{-1}$, one has
$$
\int_{\bG(\Q_p)} H({\bf s},g_p)^{-1}dg_p =
p^{-3}\left( \sum_{A\subseteq \cA} \# D^0_{A}({\bf F}_p)  
\prod_{\al\in A} \frac{p-1}{p^{s_{\al}-\kappa_{\al}+1}-1} \right),
$$
where $X=\sqcup D_A^0$ is the stratification of $X$ by locally 
closed subvarieties as in the Introduction and 
${\bf F}_p=\Z/p\Z$ is the finite field with $p$ elements. 
\end{prop}

\begin{proof}
This is Theorem 9.1 in \cite{CLT}. 
The proof proceeds as follows: 
for $p\notin S_X$ 
there is a {\em good} model $\cX$ of $X$ over $\Z_p$:
all boundary components $D_{\al}$ (and $\bG$) are defined over 
$\Z_p$ and form a strict normal crossing divisor. 
We can consider the reduction map 
$$
{\rm red}\,\,:\,\, 
X(\Q_p)=X(\Z_p)\ra X({\bf F}_p)=
\sqcup_{A\subset \cA}\, D_A^0({\bf F}_p).
$$
The main observation is that, in a neighborhood of the preimage
in $X(\Q_p)$ of the point $\tilde{x}_v\subset D_A^0({\bf F}_p)$, 
one can introduce local $p$-adic analytic coordinates 
$\{ x_{\al}\}_{\al=1,...,n}$ 
such that 
$$
H_p({\bf s},g)=\prod_{\al \in A} |x_{\al}|_p^{s_{\al}}. 
$$
Now it suffices to keep track of the change of the measure $dg_p$:
$$
dg_p= \prod_{\al\notin A} dx_{\al} \cdot
\prod_{\al \in A} |x_{\al}|_p^{k_{\al}} dx_{\al},
$$
where $dx_{\al}$ are standard Haar measures on $\Q_p$. 
The integrals obtained are elementary:
$$
\int_{{\rm red}^{-1}(\tilde{x}_p)} H_p({\bf s},g_p)^{-1}dg_p=
\prod_{\al\notin A}\int_{p \Z_p} dx_{\al} \cdot
\prod_{\al \in A}\int_{p\Z_p} 
p^{-(s_{\al}-k_{\al})v_p(x_{\al})} dx_{\al}
$$
(where $v_p(x)=\log_p(|x|_p)$ is the ordinal of $x$ at $p$). 
Summing over all $\tilde{x}_p\in X({\bf F}_p)$, we obtain the proof
(see \cite{CLT} for more details.)
\end{proof}

\no

\begin{coro}
\label{coro:0}
For all primes $p$ one has the identity
$$
\int_{\bG(\Q_p)} H_p({\bf s},g_p)^{-1}dg_p = 
\prod_{\al\in \cA} \zeta_p(s_{\al}-\kappa_{\al}+1) \cdot f_{0,p}({\bf s}),
$$
where $f_{0,p}({\bf s})$ is a 
holomorphic function in $\sT_{-1+\eps}$.
Moreover, 
there exist a $\delta>0$ and a
function $f_0({\bf s},g)$ which is 
holomorphic in $\sT_{-\delta}$ and continuous in $g\in \bG(\A)$,  
such that
$$
\zZ_0({\bf s},g)=f_0({\bf s},g)\cdot
\prod_{\al\in \cA}\zeta(s_{\al}-\kappa_{\al}+1);
$$
moreover,
$$
\lim_{{\bf s}\ra\kappa}\,\, \zZ_0({\bf s},e)\cdot 
\prod_{\al\in \cA}(s_{\al}-\kappa_{\al})=\tau(\cK_X)\neq 0,
$$
where $\tau(\cK_X)$ is the Tamagawa number defined in \cite{peyre}.
\end{coro}

\begin{proof}
Apply Corollary 9.6 in \cite{CLT}.
\end{proof}

\

\begin{nota}
\label{nota:div}
Let ${\bf a}=(a_1,a_2)\in \Q^2$ and let 
$f_{\bf a}$ be the $\Q$-rational linear form 
$$
(x,y)\mapsto a_1x+a_2y.
$$
The linear form $f_{\bf a}$ defines an 
adelic character $\eta=\eta_{\bf a}$ of $\bG(\A)$:
$$
\eta(g(x,z,y))=\psi_1(a_1x+a_2y),
$$
where again $\psi_1$ is the Tate-character of $\A/\Q$.
Write
$$
{\rm div}(\eta)=E(\eta)-\sum_{\alpha\in \cA} d_\alpha(\eta) D_\alpha
$$
for the divisor of the function $f_{\bf a}$ on the compactification $X$
(by Corollary~\ref{coro:e}, $d_{\alpha}\ge 0$ for all $\alpha\in \cA$ and $d_{\al}>0$
for at least one $\al\in \cA$).
Denote by 
$$
\begin{array}{ccc}
\cA_0(\eta) & = & \{\alpha\,|\, d_\alpha(\eta)=0\}.
\end{array}
$$

\

Let $V\subset X$ be the induced equivariant compactification of
$\bU\subset \bG$.

\begin{assu}
\label{assu:tech}
From now on we will assume that the boundary $V\setminus \bU$ 
is a strict normal crossing divisor whose components are obtained
by intersecting the boundary components of $X$ with $V$:
$$
V\setminus \bU 
=\cup_{\alpha\in \cA^V} D_{\alpha}^V=\cup_{\al\in \cA} D_{\alpha}\cap V,
$$
(with $\cA^V\subseteq \cA$).
\end{assu}

\begin{rem}
The general case can be reduced to this situation by (equivariant) 
resolution of singularities. 
\end{rem}

\

By Lemma 7.3 of \cite{CLT}, we have
$$
-K_V=\sum_{\al\in  \cA^V}\kappa^V_{\alpha} D_{\alpha}^V,
$$
with $\kappa^V_{\alpha}\le \kappa_{\al}$ (for all $\al$) and 
equality holding for $\al$ in a {\em proper} subset of $\cA$.    

\

\no
Denote by $f_a$ the $\Q$-rational linear form on the center  $\bZ$ of $\bG$
$$
z\mapsto a\cdot z.
$$ 
The linear form $f_a$ 
defines an adelic character $\psi=\psi_a$ of $\bU(\A)/\bU(\Q)$:
$$
\psi_a(g(0,z,y))=\psi_1(az).
$$
Write
$$
{\rm div}(\psi)=E(\psi)-\sum_{\alpha\in \cA^V} d_\alpha(\psi) D_\alpha
$$
for the divisor of the function  $f_a$ on $V$ and denote by  
$$
\begin{array}{ccc}
\cA_0(\psi) & = & \{\alpha\,|\, d_\alpha(\psi)=0\}.
\end{array}
$$
\end{nota}

We note that both $\cA_0(\eta) $ and $\cA_0(\psi)$ are {\em proper} subsets
of $\cA$. A precise formulation of the statement that the trivial
representation of $\bG(\A)$ is ``isolated'' in the automorphic spectrum is
the following Proposition. 

\begin{prop}
\label{prop:est-ch}
Let $\eta=\eta_{\bf a}$ and $\psi=\psi_a$ be the 
non-trivial adelic characters occuring in Proposition~\ref{prop:formm}
($a_1,a_2,a\in \frac{1}{\sn(\bfK)}\Z$). 
For any $\eps >0$ 
there exist a constant $\sc(\eps)$ and   
holomorphic bounded functions
$\phi_{\eta}(\mathbf a,\cdot)$, $\varphi_{\psi}(a,\cdot)$
on $\sT_{-1/2+\eps}$ such that, for any $\mathbf s\in \sT_{0}$, 
one has
$$
\prod_{p\notin S_{\eta}}\int_{\bG(\Q_p)} H_p(\mathbf s,g_p)^{-1}\ovl{\eta}_p(g_p)dg_p= 
\phi_{\eta}(\mathbf a,\mathbf s)\prod_{\alpha\in \mathcal A_0(\eta)}
     \zeta^{S_{\eta}}(s_\alpha-\kappa_\alpha +1);
$$ 
$$
\prod_{p\notin S_{\psi}}\int_{\bU(\Q_p)}H_p(\mathbf s,u_p)^{-1}\ovl{\psi}_p(u_p)du_p=
\varphi_{\psi}(a,\mathbf s)\prod_{\alpha\in \mathcal A_0(\psi)}
     \zeta^{S_{\psi}}(s_\alpha-\kappa_\alpha+1), 
$$ 
where  
$\zeta^S(s)=\prod_{p\notin S} (1-p^{-s})^{-1}$ 
is the incomplete Riemann zeta function.
Moreover, 
$$
\begin{array}{ccc}
|\phi_{\eta}(\mathbf a,\mathbf s)| & \le &  \sc(\eps),\\
|\varphi_{\psi}(a,\mathbf s) |     &  \le &  \sc(\eps).
\end{array}
$$
\end{prop}

\begin{proof}
The integrals can be computed as in Proposition~\ref{thm:trivial}.
They are regularized by the indicated products of (partial) zeta functions.
The remaining Euler products are 
expressions involving the number of ${\bf F}_p$-points
for boundary strata (and their intersections with ${\rm div}(\eta)$, 
resp. ${\rm div}(\psi)$).  
In particular, they are uniformly bounded on compacts in $\sT_{-1/2 +\eps}$.
For details we refer to \cite{CLT}, Proposition 5.5 
(which follows from Proposition 10.2 in {\em loc. cit.}). 
\end{proof}

\begin{coro}
In particular, 
{\em each  } term in the sums 
$\zZ_1({\bf s},g)$ and $\zZ_2({\bf s},g)$ has
a meromorphic continuation to the domain $\sT_{-1/2}$.
\end{coro}

\begin{lemm}
\label{lemm:ett}
For any $\eps>0 $ and any compact   
$\sK$ in $\sT_{-1/2+\eps}$,
there exist constants $\sc(\sK)$ and $\sn'=\sn'(\bfK)>0$, such that 
$$
\begin{array}{ccc}
|\prod_{p\in S_{\eta}} \int_{\bG(\Q_p)}
H_p({\bf s},g_p)^{-1}dg_p|& \le & \sc(\sK)\cdot (1+\|{\bf a}\|)^{\sn'}\\
|\prod_{p\in S_{\psi}} \int_{\bG(\Q_p)}
H_p({\bf s},g_p)^{-1}dg_p| & \le & \sc(\sK)\cdot (1+|a|)^{\sn'}
\end{array}
$$
for all ${\bf s}\in \sK$. 
\end{lemm}

\begin{proof}
For $p\in S_X$ we use the bound from Lemma~\ref{lemm:bad-red}. 
For $p\in S_{\eta}\setminus S_X$ 
(resp. $S_{\psi}\setminus S_X$) we apply  Proposition~\ref{thm:trivial}: 
there is a constant $\sc>0$ (depending only on $X$ and $\sK$) 
such that
$$
|\int_{\bG(\Q_p)}H_{p}({\bf s},g_p)^{-1}dg_p|\le (1+ \frac{\sc}{\sqrt{p}})
$$
for all $p$. Using the bound
$$
\prod_{p|b} (1+ \frac{\sc}{\sqrt{p}}) \le |b|^{\sn'},
$$
(for $b=a\cdot \sn(\bfK)$ and some $\sn'=\sn'(\bfK)>0$) we conclude the proof. 
\end{proof}

\begin{prop}
\label{lemm:infty}
For any $\sn>0$ and any compact $\sK\subset\sT_{-1/2+\eps}$, there exists  
a constant $\sc(\sK, \sn)$ such that, for any ${\bf s}\in \sK$,
and any $\mathbf a=(a_1,a_2)$ and $a$ as above, 
one has the estimates
$$
\begin{array}{ccl}
|\int_{\bG(\R)} H_{\infty}(\mathbf s,g_{\infty})^{-1}\ovl{\eta}_{\infty}(g_{\infty}) 
dg_{\infty}|& \leq & \sc(\sK,\sn)(1+\|{\bf a}\|)^{-\sn},
     \\
& & \\
|\int_{\bG(\A_S)} H_{S}(\mathbf s,g_{S})^{-1} 
\ovl{\omega}_{S}(g_{S}) dg_{S}| &  \leq &  \sc(\sK,\sn)(1+|\la|)^{-\sn}(1+|a|)^{\sn'},
\end{array}     
$$
where 
$\sn'=\sn'(\bfK)$ is the bound from Lemma~\ref{lemm:ett}, 
$S=S_{\psi}\cup \{\infty\}$, $\la=\la(\omega)$ 
is the eigenvalue of $\omega_{S}\in 
\cB_{S}(\varrho_{\psi})$ 
(with respect to the elliptic operator $\Delta$).
\end{prop}

\begin{proof}
We use Lemma~\ref{lemm:bad-red} and integration by parts.
For $\eta$ we apply the operator $\partial =\partial_x^2 +\partial_y^2$ 
(as in \cite{CLT})
and for $\psi$ the elliptic operator 
$\Delta=\partial_x^2+\partial_y^2+\partial_z^2$
(and use the eigenfunction property of $\omega_S$, 
or, what amounts to
the same, of $\omega\in \cB(\varrho_{\psi})$).
More precisely, 
the second integral is majorized by 
$$
|\la|^{-\sm}\cdot |\int_{\bG(\A_S)} 
\Delta^{\sm} H_S({\mathbf s},g_S)^{-1}dg_S| \cdot
\sup_{g_S\in \bG(\A_{S})} |\omega_S(g_S)|.
$$
Using the class number one property 
$$
\bG(\A)=\bG(\Q)\cdot \bG(\R)\cdot \bfK
$$
and the invariance of $\omega$ under 
$\bG(\Q)$ and $\bfK$, we obtain the estimates  
\begin{equation}
\label{eqn:12}
\sup_{g_S\in \bG(\A_{S})} |\omega_S(g_S)| \le
\sup_{g\in \bG(\A)} |\omega(g)| =
\sup_{g\in \Gamma\backslash \bG({\R})} |\omega_{\infty}(g)|.
\end{equation}
Further we have
\begin{equation}
\label{eqn:13}
\sup_{g\in \Gamma\backslash \bG({\R})} |\omega_{\infty}(g)|\ll 
|\la|^{\sm'} \cdot \|\omega\|_{\sL^2(\Gamma\backslash \bG(\R))}=
|\la|^{\sm'} \cdot \|\omega\|_{\sL^2(\bG(\Q)\backslash \bG(\A))}=|\la|^{\sm'}
\end{equation}
for some constant $\sm'$ 
(see \cite{folland}, \cite{ram}, p. 22, and \cite{sog} 
for the comparison between the $\sL^2$  
and the $\sL^{\infty}$ norms of an eigenfunction of an elliptic
operator on a compact manifold and other applications of this inequality). 
The rest of the proof follows at once 
from Lemma~\ref{lemm:bad-red} and Lemma~\ref{lemm:ett}.
(Notice that the implied constants, including $\sm'$,  
in the above inequalities depend only on the choice of $\bfK$.) 
\end{proof}

Before continuing to the proof of the main theorem, we discuss the individual
inner products
$$
\zZ({\bf s},\eta)=\langle \zZ({\bf s},\cdot ) , \eta\rangle,
$$
$$
\zZ({\bf s}, \omega^{\psi})= \langle \zZ({\bf s}, \cdot ), \omega^{\psi}\rangle.
$$ 
Let us first set
$$
\zeta_{\eta}({\bf s})=\prod_{\al\in \cA_0(\eta)} \zeta(s_{\al}-\kappa_{\al}+1),
$$
$$
\zeta_{\psi}({\bf s})=\prod_{\al\in \cA_0(\psi)} \zeta(s_{\al}-\kappa_{\al}+1).
$$
We have then the following Corollary.

\begin{coro}
\label{coro:410}
The functions 
$$
\zeta_{\eta}({\bf s})^{-1}\cdot \zZ({\bf s},\eta)
$$ 
and 
$$
\zeta_{\psi}({\bf s})^{-1}\cdot \zZ({\bf s},\omega^{\psi}),
$$
initially defined for ${\bf s}\in \sT_{\delta}$ (c.f. \ref{nota:cone}), have
an analytic continuation to the domain $\sT_{-1/2+\eps}$ (for all $\eps>0$). 
\end{coro}

\begin{proof}
We will consider the function 
$\zeta_{\psi}({\bf s})^{-1}\cdot \zZ({\bf s}, \omega^{\psi})$
and leave the first case to the reader. 
We start by observing that, for ${\bf s}\in \sT_{\delta}$, we have
\begin{equation}
\label{eqn:140}
\int_{\bG(\A)} | H({\bf s},g)|^{-1} dg < \infty
\end{equation}
(this follows from  
Proposition~\ref{prop:abs-conv} together with 
the compactness of $\bG(\Q)\ba \bG(\A)$).
Consequently,  
$$
\zZ({\bf s},\ovl{\omega}^{\psi})=\int_{\bG(\A)} H({\bf s},g)^{-1} \omega^{\psi}(g)dg,
$$
again for ${\bf s}\in \sT_{\delta}$.
Using the left-$\bfK$, and in particular, 
the left $\bfK^S$-invariance of $H$, we
have, for ${\bf s}\in \sT_{\delta}$,  
$$
\zZ({\bf s},\ovl{\omega}^{\psi})=
\int_{\bG(\A)} H({\bf s},g)^{-1}\int_{\bfK^S}\omega^{\psi}(k^Sg)dk^Sdg.
$$
Then, from Lemma~\ref{lemm:fact} (factorization), we have 
(with ${\bf s}$ in the same domain),
\begin{equation}
\label{eqn:14}
\zZ({\bf s},\ovl{\omega}^{\psi}) =\int_{\bG(\A_S)}H({\bf s}, g_S)^{-1} 
\omega^{\psi}(g_S) dg_S\cdot 
\int_{\bG(\A^S)} H({\bf s}, g^S)^{-1} f^S(g^S) dg^S,
\end{equation}
where we have set
$$
f^S(g^S):= \prod_{p\notin S} f_p(g_p)
$$
(recall that $S=S_{\psi}\cup \{\infty\}$). 
Both integrals above are convergent for ${\bf s}\in \sT_{\delta}$ by (\ref{eqn:140}). 
By Lemma~\ref{lemm:bad-red}, the first integral on the right
in (\ref{eqn:14})
is absolutely convergent
for ${\bf s}\in \sT_{-1}$. Next it follows from Proposition~\ref{prop:est-ch}
that the second integral above actually converges for ${\bf s}\in \sT_{0}$.
Moreover, we have for that integral the product expression
$$
\prod_{p\notin S_{\psi}} \int_{\bG(\Q_p)} H_p({\bf s}, g_p)^{-1} f_p(g_p)dg_p.
$$
As we have noted in Proposition~\ref{prop:est-ch}, the infinite product is  
convergent to a holomorphic function, for  ${\bf s}\in \sT_{0}$. 
Further, we then have for this infinite product the expression
$$
\varphi_{\psi}(a,{\bf s})\cdot \prod_{\al\in \cA_0(\psi)} \prod_{p\in S_{\psi}}
\zeta(s_{\al}-\kappa_{\al} +1) \cdot \zeta_{\psi}({\bf s}),
$$
for ${\bf s}\in \sT_0$. It follows, again from Proposition~\ref{prop:est-ch}, 
that 
$$
\zeta_{\psi}({\bf s})^{-1}\cdot \zZ({\bf s}, \ovl{\omega}^{\psi})
$$
can be continued holomorphically to the domain $\sT_{-1/2 +\eps}$. 
(Note that we have used 
the meromorphic continuation of $\zeta(s)$ to $\Re(s)>1/2+\eps$.)
\end{proof}

Moreover, we have the following Lemma:

\begin{lemm}
\label{lemm:413}
For any $\eps, \sn>0$ and any compact $\sK\subset \sT_{-1/2+\eps}$, there is 
a constant $\sc(\sK,\sn)$ and an integer $\sn'>0$ such that, for any ${\bf s}\in \sK$
and $a$ as above ($\psi=\psi_a$), we have
$$
|\zeta_{\psi}({\bf s})^{-1}\cdot 
\zZ({\bf s},\ovl{\omega}^{\psi})| \le \sc(\sK,\sn) 
(1+\|\lambda|)^{-\sn}(1+|a|)^{\sn'}.
$$
\end{lemm}

\begin{proof}
We have from the preceeding (proof of Corollary~\ref{coro:410})
$$
\zeta_{\psi}({\bf s})^{-1}\cdot 
\zZ({\bf s},\ovl{\omega}^{\psi}) = 
$$
$$
\varphi_{\psi}(a,{\bf s})\cdot
\int_{\bG(\A_S)} H({\bf s},g_S)^{-1}\omega(g_S)dg_S \cdot 
 \prod_{\al\in \cA_0(\psi)} \prod_{p\in S_{\psi}}
\zeta_p(s_{\al}-\kappa_{\al} +1),
$$
for ${\bf s}\in \sT_{-1/2+\eps}$. 
Our conclusion follows from 
Proposition~\ref{prop:est-ch} and Proposition~\ref{lemm:infty} and, for example, 
the elementary inequality
$$
\prod_{p|b} (1+\frac{1}{\sqrt{p}}) \le \cdot |b|^{\sn'}
$$
applied to $b=a\sn(\bfK)$
(for some $\sn'>0$, independent of $a$).
\end{proof}

\begin{thm}
\label{thm:mm}
The height zeta function 
$\zZ({\bf s})$ is holomorphic for ${\bf s}\in \sT_0$. 
Moreover, 
$$
\prod_{\al\in \cA} (s_{\al}-\kappa_{\al}) \cdot \zZ({\bf s})
$$
admits a holomorphic continuation to $\sT_{-\delta}$ (for some $\delta>0$)
and 
$$
\lim_{{\bf s}\ra \kappa}\prod_{\al\in \cA} 
(s_{\al}-\kappa_{\al}) \cdot \zZ({\bf s})=\tau(\cK_X).
$$
\end{thm}

\begin{proof}
Set 
$$
{\mathsf z}({\bf s}):=\prod_{\al \in \cA} (s_{\al}-\kappa_{\al}).
$$
We prove first that both series 
\begin{equation}
\label{eqn:16}
\sum_{\eta\neq {\bf 1}} {\mathsf z}({\bf s})\cdot \zZ({\bf s},\eta)\cdot \eta(g)
\end{equation}
and 
\begin{equation}
\label{eqn:17}
\sum_{\psi\neq 1}  \sum_{\omega^{\psi}\in \cB(\varrho_{\psi})}
{\mathsf z}({\bf s})\cdot \zZ({\bf s},\omega^{\psi})\cdot \omega^{\psi}(g) 
\end{equation}
are normally convergent for ${\bf s}$ in a compact subset of $\sT_{-1/2+\eps}$
and $g$ in a compact subset of $\bG(\A)$. We note that, by Proposition~\ref{lemm:infty}, 
the products 
$$
{\mathsf z}({\bf s})\zZ({\bf s},\eta) \,\,\, {\rm and}\,\,\, 
{\mathsf z}({\bf s})\zZ({\bf s},\omega^{\psi})
$$
are defined for ${\bf s}\in \sT_{-1/2+\eps}$.
We will prove our assertion for the second series; the proof for the first
is entirely similar.

We have a map from the set of non-trivial characters $\{ \psi\}$ of 
$\A/\Q$ to the set of subsets of $\cA$ given by
$$
\psi\mapsto \cA_0(\psi).
$$
It suffices to prove our assertion for each subseries $\zZ_A$ 
of  $\zZ_2$ corresponding to $\psi$ with $\cA_0(\psi)=A$ (for $A\subset \cA$). 
From  Lemma~\ref{lemm:413} we have a uniform majoration (for real ${\bf s}$)
$$
{\mathsf z}({\bf s})\cdot \zZ({\bf s}, \omega^{\psi})\ll
{\mathsf z}({\bf s})\cdot \zeta_{\psi}({\bf s})\cdot (1+|\la|)^{-\sn}\cdot (1+|a|)^{\sn'}.
$$
By definition, the function $\zeta_{\psi}$ is the same for all  
for $\psi$ occuring in $\zZ_A$. 
It remains then to prove the assertion for the 
series
\begin{equation}
\label{eqn:18}
\sum_{\psi} \sum_{\omega^{\psi}\in \cB(\varrho_{\psi})}
|\la|^{-\sn+\sm'}|a|^{\sn'},
\end{equation}
where we have used the estimates (\ref{eqn:12}) and (\ref{eqn:13}) (and 
the sum is over all characters $\psi$ occuring in $\zZ_A$).
We recall that 
$\la=\la(\omega^{\psi})$ is the $\Delta$-eigenvalue of $\omega^{\psi}$ and $\psi=\psi_a$. 
We also recall (Lemma~\ref{lemm:dimension}) that (with $S=S_{\psi}$) 
$$
\omega^{\psi}=j(\varphi_S\otimes e^S\otimes \varphi_n^{\psi}),
$$
for $n=0,1,2,...$, where
$\varphi_S$ varies over an orthonormal basis of $\cS(\A_S)^{\bfK_S}$. Thus our series  
(\ref{eqn:18}) is bounded from above by
$$
\sum_{a\in \Z, a\neq  0} 
\sum_{n} |\la_n|^{-\sn}|a|^{\sn'+1} \cdot \sn(\sK)^2
$$
(see Lemma~\ref{lemm:mult}).
Our claim now follows upon remarking that
$$
\la_n=(-2\pi(n+1)|a| - 4\pi^2a^2). 
$$ 

At this point we may conclude that the series (\ref{eqn:16}) and (\ref{eqn:18})
converge as stated. 
It now follows that, for ${\bf s}\in \sT_{\delta}$. 
$$
\zZ({\bf s},g)= \zZ_0({\bf s},g) + \zZ_1({\bf s},g)+\zZ_2({\bf s},g),
$$
as an equality of continuous functions on $\bG(\A)$. 
In particular, we have
\begin{equation}
\label{eqn:20}
{\mathsf z}({\bf s}) \zZ({\bf s})={\mathsf z}({\bf s})
( \zZ_0({\bf s},id) + \zZ_1({\bf s},id) + \zZ_2({\bf s},id)),
\end{equation}
again for ${\bf s}\in \sT_{\delta}$. Finally, we obtain, from (\ref{eqn:20}), 
Corollary~\ref{coro:0} and the preceding, the meromorphic continuation 
of $\zZ({\bf s})$ to the domain $\sT_{-1/2+\eps}$.

Further, since for non-trivial $\psi$ the set $\cA_0({\psi})$ 
is a {\em proper} subset of $\cA$, we also see that 
the function
$$
{\mathsf z}({\bf s}) (\zZ_1({\bf s},id)+\zZ_2({\bf s},id))
$$
vanishes for ${\bf s}=\kappa$.
Thus we have finally
$$
{\mathsf z}({\bf s}) \zZ({\bf s},id)|_{{\bf s}=\kappa} =
{\mathsf z}({\bf s}) \zZ_0({\bf s},id)|_{{\bf s}=\kappa}.
$$
Applying Corollary~\ref{coro:0} we conclude the proof.
\end{proof}

\begin{rem}
Theorem~\ref{thm:mm} implies that for each $L$ in the interior of $\La_{\rm eff}(X)$
the (one-parameter) height zeta function
$\zZ(s,\cL)$ is holomorphic for $\Re(s)>a(L)$ and admits a meromorphic continuation
to $\Re(s)>a(L)-\eps$ (for some $\eps>0$) with an isolated pole at $s=a(L)$ of
order {\em at most} $b(L)$. The proof that the order is exactly $b(L)$ and 
that the leading coefficient of $\zZ(s,\cL)$ at this pole is $c(L)\cdot \tau(\cL)$
is analogous to the proof of the corresponding statement for
height zeta functions of equivariant compactifications of additive groups
(see  Section 7 in \cite{CLT}). 
\end{rem}

\medskip

\

\section{Example: $\P^3$}
\label{sect:example}

A standard bi-equivariant compactification of the Heisenberg group $\bG$ is 
the 3-dimensional projective space $X=\P^3$. 
The boundary $D=X\setminus \bG$ 
consists of a single irreducible divisor 
(the hyperplane section). The class of this divisor 
generates the Picard group $\Pic(X)$.
The anticanonical class $-[K_X]=4 [D]$ and the cone of 
effective divisors $\Lambda_{\rm eff}(X)=\R_{\ge 0} [D]$. 
The height pairing is given by
\begin{equation}
\label{eqn:h}
H(s,g):=\prod_{p} H_p(s,g_p) \cdot H_{\infty}(s,g_{\infty}),
\end{equation}
where $g\in \bG(\A)$,
\begin{equation}
\label{eqn:p}
H_p(s,g_p)=\max \{ 1, |x|_p, |y|_p, |z|_p\}^s 
\end{equation}
and 
\begin{equation}
\label{eqn:inf}
H_{\infty}(s,g_{\infty})=(1+x^2+y^2+z^2)^{s/2}.
\end{equation}
The heights $H_p$ are invariant with respect to the 
action of $\bG(\Z_p)$ (on both sides). 
We are interested in the analytic
properties of the height zeta function
\begin{equation}
\label{eqn:zzh}
\zZ(s,g)=\sum_{\gamma=(x,z,y)\in \Q^3}H(s,\gamma g)^{-1}. 
\end{equation}

\

As above, we consider the Fourier expansion of $\zZ(s,g)$. 
Each term in this expansion will be regularized by 
an explicit Euler product of height integrals. 
We need to compute these height integrals
at good primes and estimate them at 
bad primes and at the real place.

\

\begin{lemm}
\label{lemm:trivial}
For $\Re(s)>4$, one has
$$
\int_{\bG(\A_{\rm fin})}H(s,g)^{-1}dg =\frac{\zeta(s-3)}{\zeta(s)}.
$$
\end{lemm}

\

\begin{lemm}
\label{lemm:good-eta}
For $\Re(s)>3$ and all $p\notin S_{\eta}$, one has
$$
\int_{\bG(\Q_p)}H_p(s,g_p)^{-1}\overline{\eta}_{\bf a}(g_p)dg_p = 
\zeta_p^{-1}(s).
$$
\end{lemm}

\begin{proof} 
Both Lemmas may be proved by direct computation using the 
definition of $H_p$ in (\ref{eqn:p}).
\end{proof}

\begin{lemm}
\label{lemm:good-psi}
For $\Re(s)>3$,
all $\psi=\psi_a$ and all $p\notin S_{\psi}$, 
one has
$$
\int_{\bG(\Q_p)}H_p(s,g_p)^{-1}f_p(g_p)dg_p =  
\int_{\bU(\Q_p)}H_p(s,u_p)^{-1}\overline{\psi}_{a}(u_p)du_p=\zeta_p^{-1}(s),
$$
(where $f_p$ is the local normalized spherical function).
\end{lemm}

\begin{proof} 
Direct computation. Note that the second integral
is similar to the integral in Lemma~\ref{lemm:good-eta} for the
variety $\P^2 \subset \P^3$ (the induced equivariant compactification of 
$\bU$). 
\end{proof}

\begin{lemm}
\label{lemm:bad-primes}
For all $\eps>0$, $\sn>0$ and all compacts ${\mathsf K}$ 
in the domain $\Re(s)>3+\eps$, there
exists a constant $\sc(\sn,{\mathsf K})$ 
such that, for all $s\in {\mathsf K}$ and 
all $\eta=\eta_{\bf a}$ (with ${\bf a}\in \Z^2$),
the finite product
$$
| \prod_{p\in S_{\eta}} 
\int_{\bG(\Q_p)}H_p(s,g_p)^{-1}\overline{\eta}(g_p)dg_p
\int_{\bG(\R)}H_{\infty}(s,g_{\infty})^{-1}
\overline{\eta}(g_{\infty})dg_{\infty}| 
$$
is bounded by
$$
\sc(\sn,{\mathsf K})(1+|a_1|+|a_2|)^{-\sn}.
$$
\end{lemm}

\begin{proof}
We replace $\eta$ by 1, $H_p(s,g_p)$ by
$H_p(\Re(s),g)$ and obtain 
$$
|\int_{\bG(\Q_p)}H_p(s,g_p)^{-1}dg_p|\le \frac{1}{1-p^{-\eps}}.
$$
For $a\in \Z$, we have 
$$
\prod_{p|a}(1+ p^{-\eps})\le (1+|a|)^{\sn'}
$$
(for some positive integer $\sn'$).
Using the definition of $H_{\infty}$:
$$
|\int_{\R^3}(1+x^2+y^2+z^2)^{-s/2}e^{-2\pi i(a_1x+a_2y)}dxdydz|< 
\sc(\sn,{\mathsf K})(1+|a_1|+|a_2|)^{-\sn}
$$
for all $\sn$ (this is an easy consequence of integration by parts). 
\end{proof}

\begin{lemm}
\label{lemm:bad-psiii}
For all $\eps>0$, $\sn>0$ and all compacts ${\mathsf K}$ 
in the domain $\Re(s)>3+\eps$, there
exists a constant $\sc(\sn,{\mathsf K})$ 
such that, for all $s\in {\mathsf K}$, 
all $\psi=\psi_a$ (with $a\in \Z, a\neq 0$), $S=S_{\psi}\cup \{ \infty\}$
and all $\omega_S\in \cB_S(\varrho_{\psi})$, the  expression
$$
| \int_{\bG(\A_{S})}H(s,g_S)^{-1}\overline{\omega}_S(g_S)dg_S|
$$
is bounded by
$$
\sc({\mathsf K}, \sn)|an|^{-\sn},
$$
(where the real component of 
$j_{\psi}^{-1}(\omega_S)$ is equal to  $c_n\varphi_n^{\psi}$, 
cf. Lemma~\ref{lemm:dimension}).
\end{lemm}

\begin{proof}
Let $\lambda$ be the $\Delta$-eigenvalue of $\omega$. 
Here 
$$
\Delta= \partial_x^2+\partial_y^2+\partial_z^2
$$
is an elliptic differential operator on $\bG(\Z)\ba \bG(\R)$
and 
$\partial_x$ (resp. $\partial_y,\partial_z$) 
is the invariant vector field corresponding to $g(x,0,0)$ 
(resp. $g(0,0,y)$ and $g(0,z,0)$). 
On each subspace
$$
\cH_{\psi}^{\bfK}\subset 
\sL^2(\bG(\Q)\ba \bG(\A))^{\bfK}=\sL^2(\bG(\Z)\ba\bG(\R)),
$$
we have 
$$
\partial_z \omega=(2 \pi i a)\omega
$$
(here we used the $\pi_{\psi}$ realization).
It follows that
$$
\partial_z^2\omega = -4\pi^2a^2\cdot\omega,
$$
and 
$$
\Delta\omega=(\la_n^{\psi} -4\pi^2a^2)\omega,
$$
where $\la_n^{\psi}=-2\pi(2n+1)|a|$ is the 
$\Delta_{\psi}$-eigenvalue of $\varphi_n^{\psi}$, 
the real component of $j_{\psi}^{-1}(\omega_S)$.

\

After these preparations we can assume that $s$ is real.
Using repeated integration by parts, we find the following
estimate for the above integral:
$$
\la^{-\sn}\cdot
\|\omega\|_{\sL^{\infty}}\cdot \prod_{p|a}\int_{\bG(\Q_p)}
H_{p}(s,g_p)^{-1}dg_p \cdot 
\int_{\R^3} \Delta^\sn(1+x^2+y^2+z^2)^{-s/2}dxdydz.
$$
Here we have again used the estimates (\ref{eqn:12}) and (\ref{eqn:13}). 
Continuing, we estimate the finite product of $p$-adic integrals
as in the proof of Lemma~\ref{lemm:bad-primes}. Finally, we find 
from Lemma~\ref{lemm:bad-red} that the integral
$$
\int_{\R^3} \Delta^\sn(1+x^2+y^2+z^2)^{-s/2}dxdydz
$$ 
is convergent for $s\in \sK$ and further is bounded 
on the same region. 
\end{proof}

\begin{prop}
\label{prop:zetah}
The height zeta function $\zZ(s)$ defined 
in Equation~(\ref{eqn:zzh})
\begin{itemize}
\item is holomorphic for $\Re(s)>4$;
\item admits a meromorphic continuation to $\Re(s)>3+\eps
$ (for any $\eps>0$) and 
\item has a simple pole in this domain at $s=4$
with residue
$$
\zeta(4)^{-1}\int_{\R^3}(1+x^2+y^2+z^2)^{-2}dxdydz=\frac{\pi^2}{\zeta(4)}.
$$
\end{itemize}
\end{prop}

\begin{proof}
Using the estimates of Lemma~\ref{lemm:bad-psiii}, and 
(\ref{eqn:12}) and (\ref{eqn:13}), we see, as in the proof of Theorem\ref{thm:mm}, 
that the series for $\zZ_2(s,g)$ in Proposition~\ref{prop:formm}
is normally convergent for $s$ in a compact set $\sK$ of $\Re(s)>3$ and
for $g\in \bG(\A)$. It now follows (as in the proof of Theorem~\ref{thm:mm})
that, for $\Re(s)>4$, 
$$
\zZ(s,g)=\zZ_0(s,g)+\zZ_1(s,g)+\zZ_2(s,g),
$$
as an equality of continuous functions on $\bG(\A)$.
In particular,
\begin{equation}
\label{eqn:21}
 \zZ(s,id)=\zZ_0(s,id)+\zZ_1(s,id)+\zZ_2(s,id),
\end{equation}
again for $\Re(s)>4$. 
We obtain then, as in the proof of Theorem~\ref{thm:mm}, 
the meromorphic continuation of $\zZ(s)$ to $\Re(s)>3+\eps$
(see esp. Lemma~\ref{lemm:trivial} for $\zZ_0$). 
Finally, 
$$
(s-4)\zZ(s)|_{s=4}=(s-4)\zZ_0(s,id)|_{s=4}=\frac{\pi^2}{\zeta(4)}.
$$
\end{proof}

\end{document}